\documentclass[11pt]{article}

\usepackage{amsmath,amssymb,latexsym}

\title{Kreisel's counter-example to full abstraction of the \\ 
  set-theoretical model of G\"odel's system $T$}

\author{Mart\'\i n Escard\'o}

\date{23rd July 2007}

\newcommand{\N}{\mathbb{N}}
\newcommand{\Ni}{\mathbb{N}_\infty}
\newcommand{\test}{\operatorname{test}}

\begin{document}

\maketitle

\paragraph{Introduction.}
I have written this note because it doesn't seem to be well known that
the set-theoretical model of G\"odel's system $T$ fails to be fully
abstract. By the context lemma, full abstraction is equivalent to the
statement that the substructure of definable elements is extensional;
that is, that any two definable functionals that agree on definable
arguments must agree on all arguments. Kreisel gave a counter-example
to this, reported in~\cite[page 581, Exercise 1]{barendregt} by
Barendregt. A ground-valued definable functional that is constant on
definable arguments but non-constant on arbitrary arguments is
exhibited.

The surprise in the construction is that exhaustive search over the
one-point compactification of the discrete natural numbers is
definable in system~$T$. Once one knows this, it is straightforward to
produce the counter-example, using the well-known fact that definable
functionals $2^\N \to 2$ are continuous with respect to the discrete
topology on $2$ and the product topology on $2^\N$. The connection
with exhaustive search is in fact made explicit in Exercise~2
of~\cite[page 581]{barendregt}. We propose a slight improvement of the
suggested solution, and use this to formulate Kreisel's
counter-example.

\paragraph{System $T$.} We take system $T$ to be the simply typed
lambda-calculus with base types for natural numbers (including zero,
successor and higher-type primitive recursion) and booleans (including
if-then-else).

\paragraph{The set-theoretical model.}
The set theoretical model interprets the type of natural numbers as
the set $\N$ of natural numbers and the type of booleans as the set
$2=\{0,1\}$, with the convention that $0$ is false and $1$ is true.
Function types are interpreted as exponentials in the category of sets
(=sets of all functions).

\paragraph{The one-point compactification of $\N$.} 
Let \[ \Ni = \{ x \in 2^\N \mid \forall i\le j.x_i \le x_j \} = \{
\bar{n} \mid n \in \N \} \cup \{ \infty \}, \] where we write $\bar{n}
= 0^n1^\omega$ and $\infty = 0^\omega$. The last equation requires
classical logic but we don't worry about this (see below). The
argument given below uses the fact that definable functions $2^\N \to
2$ are continuous with respect to the discrete topology on $2$ and the
product topology on $2^\N$ (Cantor space), and the fact that $\infty =
\lim_n \bar{n}$ in this topology. Notice that the pointwise order on
$\Ni$ restricts to the natural order on $\N$ and has $\infty$ as a top
element.  For $\alpha \in 2^\N$, denote by $\alpha+1$ the sequence
$0\alpha$.  Then $\bar{n}+1 = \overline{n+1}$ and $\infty + 1 =
\infty$.

\paragraph{Exhaustive searchability of $\Ni$.} 
Define $\varepsilon_{\Ni} \colon 2^{2^\N} \to 2^\N$ by
\[
\text{$\varepsilon_{\Ni}(p)(i) = 1$ iff $p(\bar{n})$ for some $n \le i$.}
\] 
Because the bounded existential quantification can be reduced to
primitive recursion, this is $T$-definable. By construction,
\[
\text{$\varepsilon_{\Ni}(p) = \overline{\mu n.p(\bar{n})}$ if
  $p(\bar{n})$ for some $n \in \N$,}
\]
and
\[
\text{$\varepsilon_{\Ni}(p) = \infty$ if there isn't $n \in \N$ such
  that $p(\bar{n})$.}
\]
Or, combining these two statements,
\[
\varepsilon_{\Ni}(p) = \inf \{ x \in \Ni \mid p(x) \}
\]
because the infimum of the empty set is always the top element.
Hence the image of $\varepsilon_{\Ni}$ is $\Ni$, and
\[
\text{$p(\varepsilon_{\Ni}(p))$ iff $p(x)$ for some some $x
  \in \Ni$.}
\]
In the terminology of~\cite{escardo:lics2007}, the set $\Ni$ is
searchable.  Notice that $\varepsilon_{\Ni}$ implements search even for
\emph{discontinuous}~$p$. 

\pagebreak[4]

\paragraph{Kreisel's counter-example.} Define $\test \colon
2^{2^\N} \to 2$ by
\[
\test(f) = f\left[\inf x \in \Ni.f(x+1) = f(\infty)\right] = f(\infty) \implies
f(\bar{0}) = f(\infty).
\]
By the above, this is $T$-definable. If $f \colon 2^\N \to 2$ is
continuous, then $\test(f)$ holds, because there is $n \in \N$ such
that $f(\bar{n})=f(\infty)$ as $\lim_n \bar{n} = \infty$. Hence the
functional $\test$ is constant on $T$-definable arguments. However,
$\test(f)$ fails for $f$ defined by
\[
f(\alpha) =
\begin{cases}
  1 & \text{if $\alpha = \infty$,} \\
  0 & \text{otherwise,}
\end{cases}
\]
because in this case $[\inf x \in \Ni.f(x+1) = f(\infty)] = \infty$
and hence 
\[
\text{$f\left[\inf x \in \Ni.f(x+1) = f(\infty)\right] = f(\infty)$
  but $f(\bar{0}) \ne f(\infty)$.}
\]

\paragraph{A fully abstract model.} Kleene--Kreisel functionals
(exponentials in $k$-spaces starting from discrete natural numbers and
booleans) are an example, using the fact that the Kleene--Kreisel
density theorem gives dense sequences which are actually definable in
system $T$.

\paragraph{Classical logic.} The above argument relies on classical
logic (in a weak form: a sequence is either $\infty$ or not). This is
necessarily the case. From a fully abstract model, taking presheaves
one gets a topos (model of bounded set theory based on intuitionistic
logic) that gives an equivalent fully abstract model!

\paragraph{Acknowledgements.} 
When I posed the question, Gordon Plotkin answered that Barendregt
knew a counter-example in the late 1970's, and Alex Simpson later told
me where to find it.

\bibliographystyle{plain}
\bibliography{references}

\end{document}